\documentclass[11pt]{article}
\usepackage{amsmath,amssymb}

\newtheorem{theorem}{Theorem}
\newtheorem{corollary}[theorem]{Corollary}

\begin{document}

\title{Convexity of geodesic-length functions: a reprise}         % Enter your title between curly braces
\author{Scott A. Wolpert}        % Enter your name between curly braces
\date{May 24, 2004}          % Enter your date or \today between curly braces
\maketitle

\begin{abstract}
New results on the convexity of geodesic-length functions on Teichm\"{u}ller space are presented.  A formula for the Hessian of geodesic-length is presented. New bounds for the gradient and Hessian of geodesic-length are described.  A relationship of geodesic-length functions to Weil-Petersson distance is described.  Applications to the behavior of Weil-Petersson geodesics are discussed.
\end{abstract}

\section{Introduction}
  % replace this with a key based on an author name and refer to it in \firstpage and \lastpage
                        % above, to get page range correct

In this research brief we describe a new approach to the work \cite{Wlnielsen}(esp. Secs. 3 and 4), as well as new results and applications of the convexity of geodesic-length functions on the Teichm\"{u}ller space $\mathcal T$.  Our overall goal is to obtain an improved understanding of the convexity behavior of geodesic-length functions along Weil-Petersson (WP) geodesics. Applications are presented in detail for the $CAT(0)$ geometry of the augmented Teichm\"{u}ller space.  A complete treatment of results is in preparation \cite{Wlconvexity}.  Convexity of geodesic-length functions has found application for the convexity of Teichm\"{u}ller space \cite{Brkwpvs,Brkwp,SerDia,Kerck,Kerckmin,McM,SS2,SS1,Wlnielsen}, for the convexity of the WP metric completion \cite{DW2,MW,Wlcomp,Yam}, for the study of harmonic maps into Teichm\"{u}ller space \cite{DW1,Yam1,Yam}, and for the action of the mapping class group \cite{DW2,MW}.   We consider marked Riemann surfaces $R$ with complete hyperbolic metrics possibly with cusps and consider the lengths of simple closed geodesics.  The length of the unique geodesic in a prescribed free homotopy class provides a function on the  Teichm\"{u}ller space.  Specifically  for $\sigma$ a simple closed curve on $R$, let $\ell_{\sigma}(R)$ denote the length of the 
geodesic homotopic to $\sigma$; more generally for $\mu$ a measured geodesic lamination \cite{Bonmgl}, let $\ell_{\mu}(R)$ denote the total-length of the geodesic lamination on $R$.

A closed geodesic $\sigma$ on  $R$ determines a cyclic cover of $R$ by a geometric cylinder $\mathcal C$.  For $\ell$ the length of $\sigma$ the geometric cylinder is represented as $\mathbb H/<t\rightarrow e^\ell t>$ for $\mathbb H$ the upper half-plane with coordinate $t$; for $w=\exp(2\pi i\frac{\log t}{\ell})$ the cylinder is further represented as the concentric annulus $\{e^{-\frac{2\pi^2}{\ell}}<|w|<1\}$ in the plane.  We discovered in \cite{Wlnielsen}(Sec. 4) that the potential operator for the Beltrami equation on $\mathcal C$ is diagonalized by the $\mathbf S^1$ rotation action of the cylinder and that the potential equation can be solved term-by-term for the corresponding Fourier expansions.  The special properties for the potential theory generalize the properties for the function theory of the cylinder.  For instance holomorphic differentials on $R$, lifted to $\mathcal C$, admit Laurent (Fourier) expansions.  The WP dual of the Hessian of $\ell_{\sigma}$, a quadratic form for holomorphic quadratic differentials on $R$, has Hermitian and complex-bilinear components {\em diagonalized} by the terms of the corresponding Laurent expansions \cite{Wlnielsen}(see Lemmas 4.2 and 4.4.)  We further found that the contribution for a single Laurent term is a positive definite form.  At this time, we have simplified the considerations of the Hessian and are now able to effect a straightforward comparison to the Petersson pairing for holomorphic quadratic differentials \cite{Wlconvexity}. The simplified considerations provide the basis for an improved understanding of the Hessian and of convexity.  In the following paragraphs we outline the approach and results.  We close the discussion by providing two applications complete with proofs.  

\section{The Hessian of geodesic-length}

We introduce for $\mu$ a measured geodesic lamination  a natural function $\mathbb P_{\mu}$ on $R$.  We begin with the geometry of the space of complete geodesics on the hyperbolic plane.  For $\mathbb H$ the upper half plane with boundary $\Check{\mathbb R}=\mathbb R\cup\{\infty\}$,  the space of complete geodesics on $\mathbb H$ is given as $\mathcal G=\Check{\mathbb R}\times\Check{\mathbb R}\setminus\{diagonal\}/\{interchange\}$.  A point $p$ of $\mathbb H$ is at finite distance $d(p,\sigma)$ to a complete geodesic $\sigma$ and so $e^{-2d(p,\sigma)}$ defines a {\em Gaussian} on $\mathcal G$.  The natural area measure on $\mathcal G$ is $\omega=(a-b)^{-2}da\,db$ in terms of the {\em endpoint coordinates} $(a,b)/\sim$.  The measure $e^{-2d(p,\sigma)}\omega$ is finite for $\mathcal G$.  Finiteness is noted as follows.  A point $z$ of $\mathbb H$,  its conjugate $\bar z$, and the boundary points $(a,b)$ have cross ratio $cr(z,a,b)=\frac{(a-b)\,\Im z}{|z-a||z-b|}$.  
The simple inequality $cr^2(z,a,b)\ge e^{-2d(z,\stackrel{\frown}{ab})}$ is established by considering the point triple $(i,a,-a)$.  Finiteness of the measure now follows from the inequality $e^{-2d(i,\stackrel{\frown}{ab})}\omega\le (1+a^2)^{-1}(1+b^2)^{-1}da\,db$ for the point triple $(i,a,b)$.  A measured geodesic lamination $\mu$ on $R$ naturally lifts to the upper half plane; the lift determines a measure $d\mu$ on the space $\mathcal G$ of complete geodesics.  A compact arc on $\mathbb H$ transverse to the lamination determines a set $\tau$ of intersected leaves, a subset  of $\mathcal G$, with $\mu(\tau)$ the measure of the subset.  For $R$ represented as the quotient $\mathbb H/\Gamma$  the integral
\begin{equation}
\label{Pmu}
\mathbb P_{\mu}(p)=\int_{\mathcal G}e^{-2d(p,\sigma)}d\mu(\sigma)
\end{equation}
defines a $\Gamma$-invariant function on $\mathbb H$, {\em the mean-squared inverse exponential-distance of $p$ to} $\mu$.
Finiteness of the integral is established by comparing $d\mu$ to $\omega$.   The construction for $\mathbb P_{\mu}$ is motivated by the construction for the classical Petersson series representing the differential $d\ell_{\sigma}$ of the geodesic-length on $\mathcal T$ \cite{Gardtheta,Gardmeas}.  The reader can check that  $\mu\rightarrow\mathbb P_{\mu}$ is a continuous mapping from the space of measured geodesic laminations to the space of continuous functions on $R$.  The central role of $\mathbb P_{\mu}$ in studying geodesic-length functions and the total-length of laminations is discussed and demonstrated below.     

From Kodaira-Spencer deformation theory the infinitesimal deformations of $R$ are represented by the Beltrami differentials $\mathcal H(R)$ harmonic with respect to the 
hyperbolic metric \cite{Ahsome}.  A harmonic Beltrami differential is the symmetric tensor given as $\overline{\varphi}(ds^2)^{-1}$ for $\varphi$ a holomorphic quadratic differential with at most simple poles at the cusps and $ds^2$ the hyperbolic metric tensor.  At $R$ the differential on $\mathcal T$ of the geodesic-length of $\ell_{\sigma}$ is bounded for $\nu\in\mathcal H(R)$ as
$$
|d\ell_{\sigma}(\nu)|\le \frac{8}{\pi}\int_R|\nu|\,\mathbb P_{\sigma}dA
$$
for $dA$ the hyperbolic area element and from applying the inequality 
$|(\frac{\Im z}{\bar z})^2|\le 4 e^{-2d(z,\stackrel{\frown}{0\infty})}$ and the formula of F. Gardiner \cite{Gardtheta}.  By taking limits the integral bound is generalized to the total-length of laminations.  Ahlfors  noted \cite{Ahsome} for second-order deformations defined by harmonic Beltrami differentials that the WP Levi-Civita connection is  Euclidean to zeroth order in the following sense.  A $\Gamma$-invariant Beltrami differential $\nu$ on $\mathbb H$ determines a one-parameter family as follows.  For the complex parameter $\epsilon$ small there is a suitable self-homeomorphism $f^{\epsilon}$ of $\mathbb H$ satisfying $f^{\epsilon}_{\bar z}=\epsilon\nu f^{\epsilon}_z$.  The homeomorphism $f^{\epsilon}$ serves to compare the quotients $\mathbb H/\Gamma$ and $\mathbb H/f^{\epsilon}\circ \Gamma\circ (f^{\epsilon})^{-1}$.  For a basis of harmonic Beltrami differentials $\nu_1,\dots,\nu_n$ and small complex parameters $\epsilon_*$ and $\nu(\epsilon)=\sum_j\epsilon_j\nu_j$ the association $(\epsilon_1,\dots,\epsilon_n)$ {\em to} $\mathbb H/f^{\nu(\epsilon)}\circ\Gamma\circ(f^{\nu(\epsilon)})^{-1}$ in effect provides a local coordinate for $\mathcal T$.  Ahlfors found for a basis of harmonic Beltrami differentials that the local coordinates for $\mathcal T$ are normal: the first derivatives of the WP metric tensor vanish at the origin \cite{Ahsome}.  The observation is used in the calculation of the WP Riemannian Hessian $\Ddot{\ell}_{\mu}(\nu,\nu)$.  

Our analysis of the Hessian consists of three considerations.  We consider the metric cover of the cylinder $\mathcal C$ by an infinite horizontal strip $\mathcal S$ in $\mathbb C$ with the $\mathbf S^1$ rotation action of the cylinder lifting to an $\mathbb R$ action by Euclidean translations of the strip.  We purposefully normalize the covering $\mathcal S$ so that a Euclidean horizontal translation by $\delta$ is a hyperbolic isometry with translation length $\delta$.
First, we consider the formula for the variation of the translation length $\ell$ of the covering of $\mathcal C$.  For $z$ the complex coordinate for the strip and $f^{\epsilon}$ the suitable self-homeomorphism of $\mathcal S$ the translation equivariance provides that $\ell^{\epsilon}=f^{\epsilon}(z+\ell)-f^{\epsilon}(z)$.  We find for $\nu$ a harmonic Beltrami differential defining a deformation and $\mathcal F$ a fundamental domain for the metric covering of $\mathcal S$ to $\mathcal C$ the first variation
$$
\dot\ell=\frac{1}{\pi}Re\int_{\mathcal F}\frac{d}{d\epsilon}f^{\epsilon}_{\bar z}\,idzd\bar{z}
=\frac{1}{\pi}Re\int_{\mathcal F}\nu \,idzd\bar{z}
$$
and the second variation
$$
\Ddot\ell=\frac{1}{\pi}Re\int_{\mathcal F}\frac{d^2}{d\epsilon^2}f^{\epsilon}_{\bar z}\,idzd\bar{z}=\frac{2}{\pi}Re\int_{\mathcal F}\nu f_z\,idzd\bar{z}
$$
for $f$ a suitable solution of the potential equation $f_{\bar z}=\nu$.  The second variation formula should be compared to the considerably more involved formula of Theorem 3.2 of \cite{Wlnielsen}.  Second, we consider the Fourier expansion of $\nu$ on $\mathcal S$ relative to the translation group of the covering to $\mathcal C$.  From Corollary 2.5 and formulas (4.1) of \cite{Wlnielsen} the potential equation $f_{\bar z}=\nu$ admits a term-by-term solution relative to the Fourier expansion of $\nu$. In particular for the Beltrami differential with series expansion 
$$
\nu=-4\sin^2\Im z\ \overline{\sum a_ne^{\epsilon nz}},\ \epsilon =\frac{2\pi i}{\ell},
$$ 
we find that
$$
f_z=2\bigl(e^z\Re \sum a_n\frac{e^{\epsilon nz-1}}{\epsilon n-1}- e^{-z}\Re \sum a_n\frac{e^{\epsilon nz+1}}{\epsilon n+1}\bigr).
$$
The quantity $f_z$ is a linear form in the Fourier expansion of $\nu$. The expansion enables calculation of the above variation integral term-by-term  and the calculation is a special feature for harmonic Beltrami differentials.  Third, we simplify the resulting term-by-term expressions to obtain an exact formula in terms of the operator
$$
A[\varphi]=\zeta^{-1}\int^\zeta t^2\varphi\,dt
$$
for quadratic differentials $\varphi$ invariant by $t\rightarrow e^{\ell}t$ on $\mathbb H$ with coordinate $t$, and the Hermitian form
$$
Q(\beta,\delta)=\int_{1<|t|<e^{\ell}}\beta\bar\delta\,(Im\, t)^2\,\frac{i}{2}dtd\bar t.
$$
In \cite{Wlnielsen}(Thrm. 2.4) we found that $A[\varphi]$ is associated to the Eichler integral of $\varphi$.
The overall resulting final formula
\begin{equation}
\Ddot\ell=\frac{32}{\pi}Q(A,A)-\frac{16}{\pi}Q(A,\bar A)
\end{equation}
is the replacement for the intricate formulas of Lemmas 4.2 and 4.4 of \cite{Wlnielsen}.  The formula can be compared to the formula of Gardiner for the first variation \cite{Gardtheta}.  Bounds for the Hessian of $\ell_{\sigma}$ in terms of $\mathbb P_{\sigma}$ and the Petersson product can be derived by comparing the two Hermitian forms
$$
Q(A,A)
\quad\mbox{and}\quad
\int_{1<|t|<e^\ell}|\varphi|^2(Im\, t)^4\,\frac{i}{2} \frac{dtd\bar t}{|t|^2}
$$
where $(Im\,t)^{-1}|t|$ is comparable to the exponential-distance of $t$ to the imaginary axis.  The bounds are straightforward since the Hermitian forms are diagonalized by the Fourier expansion of $\varphi$.

\section{Convexity results}

We find that for the total-length $\ell_{\mu}$ of a measured geodesic lamination its complex Hessian on $\mathcal T$, a Hermitian form on $\mathcal H(R)$, is bounded in terms of the integral pairing with factor $\mathbb P_{\mu}$ and hyperbolic area element
$$
\int_R\nu\overline{\rho}\,\mathbb P_{\mu}dA\le \frac{3\pi}{16}\partial\overline{\partial}\ell_{\mu}(\nu,\rho)\le 16\int_R\nu\overline{\rho}\,\mathbb P_{\mu}dA
$$
for $\nu,\rho\in\mathcal H(R)$.  Since $\int_R\nu\overline{\rho}\,dA$ is the WP pairing $\bigl<\nu ,\rho\bigr>_{WP}$, we have the following comparison of Hermitian forms 
$$
\bigl<\ ,\ \mathbb P_{\mu}\bigr>_{WP}\le \frac{3\pi}{16}\partial\overline{\partial}\ell_{\mu}\le 16\bigl<\ ,\ \mathbb P_{\mu}\bigr>_{WP}.
$$
The strict convexity of geodesic-length functions and of the total-length of laminations is an immediate consequence of the positivity of $\mathbb P_{\mu}$.  
We find further consequences of our calculations and considerations of $\mathbb P_{\mu}$.  The first and second derivatives of total-lengths $\ell_{\lambda}\,,\ell_{\mu}$ actually satisfy general comparisons
\begin{equation}
\label{firstbound}
|d\ell_{\lambda}(\nu)d\ell_{\mu}(\nu)| < \ell_{\lambda}\Ddot{\ell}_{\mu}(\nu,\nu)+\ell_{\mu}\Ddot{\ell}_{\lambda}(\nu,\nu)
\end{equation}
and
\begin{equation}
\label{firstcbound}
4|\partial\ell_{\lambda}(\nu)\overline{\partial}\ell_{\mu}(\nu)|<
\ell_{\lambda}\partial\overline{\partial}\ell_{\mu}(\nu,\nu)+\ell_{\mu}\partial\overline{\partial}\ell_{\lambda}(\nu,\nu).
\end{equation}
The complex Hessian and WP Riemannian Hessian of a total-length $\ell_{\mu}$ also satisfy a general comparison
$$
\partial\overline{\partial}\ell_{\mu}\le \,\Ddot{\ell}_{\mu} \le 3\,\partial\overline{\partial}\ell_{\mu}.
$$

A basic consequence of the formulas is the observation that the first and second derivatives of a geodesic-length $\ell_\sigma$ are bounded in terms of the supremum norm of $\mathbb P_{\sigma}$ on $R$.  The magnitude of $\mathbb P_{\sigma}$ can in turn be analyzed in terms of the {\em thick-thin} decomposition of the surface \cite[II, Sec. 2]{Wlspeclim}.  For $\sigma$ a simple closed geodesic a suitable decomposition of $R$ has three regions: {\em i)} {\em thick} ; {\em ii)} cusps and {\em thin} collars not intersecting $\sigma$; and {\em iii)} {\em thin} collars which $\sigma$ crosses.  For the first region since $e^{-2d(p,\sigma)}$ satisfies a mean value estimate and the injectivity radius is uniformly bounded below the supremum of $\mathbb P_\sigma$ is bounded by the $L^1$-norm $\|\mathbb P_{\sigma}\|=\frac43\ell_{\sigma}$.  For the second region the distance to $\sigma$ is at least the distance $\delta$  to the region boundary and the supremum can be bounded using the general inequality $e^{\delta}\rho>c$ bounding the exponential-distance and the injectivity radius for a collar or cusp \cite[II, Lem. 2.1]{Wlspeclim}.  For the third region the supremum of $\mathbb P_{\sigma}$ is bounded in terms of the reciprocal injectivity radius, which from the general inequality is bounded by $e^{\ell_{\sigma}/2}$ since $\sigma$ crosses the collar.  

We accordingly find in complete generality that there exists constants $c_*,c_{**}$ independent of $R$ such that the gradient of the geodesic-length of a simple curve is bounded in terms of the geodesic-length itself
\begin{equation}
\label{gradbd}
\bigl<grad\,\ell_\sigma,grad\,\ell_\sigma\bigr>_{WP} \le c_*(\ell_{\sigma}+\ell_{\sigma}^2e^{\ell_{\sigma}/2})
\end{equation}
and for the relative systole $sys_{rel}(R)$, the least (closed) geodesic-length for $R$, that
\begin{equation}
\label{hessbd}
c_{**}(sys_{rel}(R))^{4\dim_{\mathbb C}\mathcal T}\bigl<\ ,\ \bigr>_{WP}\le \partial\overline{\partial}\ell_{\sigma}\le c_*(1+\ell_{\sigma}e^{\ell_{\sigma}/2})\bigl<\ ,\ \bigr>_{WP}.
\end{equation}
In brief the first and second derivatives of a geodesic-length relative to the WP metric are universally bounded in terms of the geodesic-length and the relative systole. The bound (\ref{gradbd}) can be compared to the familiar universal bound $\|d\ell_{\sigma}\|_T\le2\ell_{\sigma}$ for the differential relative to the Teichm\"{u}ller metric \cite{Gardtheta}. The {\em degeneration} of $\mathbb P_{\sigma}$ can be further analyzed in terms of the {\em thick-thin} decomposition \cite[II, Sec. 2]{Wlspeclim}.  

For the study of measured geodesic laminations and applications of geodesic-lengths it is desirable to have bounds (dependent on $R$) proportional to the geodesic-length.  We find for compact subsets of the moduli space of Riemann surfaces that there are uniform bounds.  In particular we have the following.
\begin{theorem}
\label{convex}
Given $\mathcal T$, there are functions $c_1$ and $c_2$ such that for a simple curve $\sigma$ 
$$
c_1(sys_{rel}(R))\,\ell_{\sigma}\le \mathbb P_{\sigma}\le c_2(sys_{rel}(R))\,\ell_{\sigma}
$$
with $c_1(s)$ an increasing function vanishing at the origin and $c_2(s)$ a decreasing function.  For the total-length of a geodesic lamination $\mu$ $$c_1(sys_{rel}(R))\,\ell_{\mu}\bigl<\ ,\ \bigr>_{WP}\le\partial\overline{\partial}\ell_{\mu}
\le c_2(sys_{rel}(R))\,\ell_{\mu}\bigl<\ ,\ \bigr>_{WP}.
$$
\end{theorem}

The first-derivative second-derivative comparison inequalities (\ref{firstbound}), (\ref{firstcbound}) provide for new convexity results.
\begin{theorem}
For the closed curves $\alpha_1,\dots,\alpha_n$ their geodesic-length sum $\ell_{\alpha_1}+\cdots+\ell_{\alpha_n}$ satisfies:   $(\ell_{\alpha_1}+\cdots+\ell_{\alpha_n})^{1/2}$ is strictly convex along WP geodesics, $\log (\ell_{\alpha_1}+\cdots+\ell_{\alpha_n})$ is strictly plurisubharmonic, and $(\ell_{\alpha_1}+\cdots+\ell_{\alpha_n})^{-1}$ is strictly plurisuperharmonic.
\end{theorem}

C. McMullen found and used that an $\ell_{\alpha}^{-1}$ has complex Hessian uniformly bounded relative to the Teichm\"{u}ller metric \cite[Thm. 3.1]{McM}.  The present result offers an elaboration: a sum $(\ell_{\alpha_1}+\cdots+\ell_{\alpha_n})^{-1}$ is plurisuperharmonic with complex Hessian bounded as 
\begin{equation}
-\partial\overline{\partial}((\ell_{\alpha_1}+\cdots+\ell_{\alpha_n})^{-1})<(2\partial\overline{\partial}(\ell_{\alpha_1}+\cdots+\ell_{\alpha_n}))(\ell_{\alpha_1}+\cdots+\ell_{\alpha_n})^{-2}.
\end{equation}  

\section{The $CAT(0)$ geometry of the augmented Teichm\"{u}ller space}

Applications of geodesic-length convexity are provided by considering the augmented Teichm\"{u}ller space $\overline{\mathcal T}$ with the completion of the WP metric \cite{Abdegn, Bersdeg,Msext}.  $\overline{\mathcal T}$ is the space of marked possibly noded Riemann surfaces; $\overline{\mathcal T}$ is a non locally compact space \cite{Abdegn, Bersdeg}.  $\overline{\mathcal T}$ is a $CAT(0)$ metric space \cite{DW2, MW, Wlcomp, Yam}.  The geometry of $CAT(0)$ spaces is developed in detail in Bridson-Haefliger \cite{BH}. For a metric space a {\em geodesic triangle} is 
prescribed by a triple of points and a triple of joining length-minimizing curves. A characterization of curvature 
for metric spaces is provided in terms of distance-comparisons to geodesic triangles in constant curvature spaces.  For a $CAT(0)$ space the distance and angle measurements for a triangle are bounded by the corresponding measurements for a Euclidean 
triangle with the corresponding edge-lengths \cite[Chap. II.1, Prop. 1.7]{BH}.

$\overline{\mathcal T}$ with the completion of the WP metric is  a {\em stratified} unique geodesic space with the strata intrinsically 
characterized by the metric geometry \cite[Thm. 13]{Wlcomp}.  The stratum containing a given point is the union of all open length-minimizing segments containing the point.
To characterize the strata structure in-the-large consider a reference topological surface $F$ for the marking and $C(F)$, the partially ordered set {\em the complex of curves}.  A $k$-simplex of $C(F)$ consists of $k+1$ distinct nontrivial
free homotopy classes of nonperipheral mutually disjoint simple closed curves. Consider $\Lambda$ the
natural labeling-function from $\overline{\mathcal T}$ to $C(F)\cup\{\emptyset\}$.  For a
marked noded Riemann surface $(R,f)$ with $f:F\rightarrow R$, the labeling
$\Lambda((R,f))$ is the simplex of free homotopy classes on $F$ mapped to the
nodes on $R$. The level sets of $\Lambda$ are the strata of $\overline{\mathcal T}$ \cite{Abdegn, Bersdeg}.  The strata of $\overline{\mathcal T}$ are lower-dimensional Teichm\"{u}ller spaces; each stratum  with its natural WP metric isometrically embeds into the completion $\overline{\mathcal T}$ \cite{Msext}.  The
unique WP geodesic $\widehat{pq}$ connecting $p,q\in \overline{\mathcal T}$ is contained in
the closure of the stratum with label $\Lambda(p)\cap\Lambda(q)$ (see \cite[Thm. 13]{Wlcomp}). 
The open segment $\widehat{pq}-\{p,q\}$ is a solution of the
WP geodesic differential equation on the stratum with label
$\Lambda(p)\cap\Lambda(q)$.  It follows from Theorem \ref{convex} that a geodesic-length function finite on $\widehat{pq}$ is necessarily strictly convex and on the open segment differentiable.  

A complete, convex subset $\mathcal C$ of a $CAT(0)$ space is the base for an
{\em orthogonal projection}, \cite[Chap. II.2, Prop. 2.4]{BH}. For a general point $p$ there is a unique point, {\em the projection of} $p$, on $\mathcal C$ such that the connecting geodesic realizes the distance to $\mathcal C$.  The projection is a retraction that does not increase distance.
The distance $d_{\mathcal C}$ to $\mathcal C$ is a convex function satisfying
$|d_{\mathcal C}(p)-d_{\mathcal C}(q)|\le d(p,q)$, \cite[Chap. II.2, Prop.
2.5]{BH}.  Examples of complete, convex sets $\mathcal C$ are: points,
complete geodesics, and fixed-point sets of isometry groups. In the case of $\overline{\mathcal T}$ since geodesics coincide at most at endpoints, the fibers of a
projection are filled out by the geodesics realizing distance between their points
and the base.   In
the case of $\overline{\mathcal T}$ the closure of each
individual stratum is complete and convex, thus the base of a projection \cite[Thm. 13]{Wlcomp}.
For simple disjoint closed curves the relation of the quantity $\ell^{1/2}=(\ell_{\alpha_1}+\cdots+\ell_{\alpha_n})^{1/2}$ to a stratum of $\overline{\mathcal T}$ was considered in \cite[Cor. 21]{Wlcomp}.  The expansion of the WP metric about a stratum \cite[Cor. 4]{Wlcomp} enabled us to give an expansion for the distance to a stratum.  The expansion combines with the comparison inequality (\ref{firstbound}) to provide an inequality for distance in-the-large.  In the following the quantity $\ell^{1/2}$ serves as a {\em Busemann function} for the stratum of vanishing.
\begin{theorem}
\label{dist}
For closed curves $\alpha_1,\dots,\alpha_n$ represented by simple disjoint distinct free homotopy classes, let $\mathcal S$ be the closed stratum of $\overline{\mathcal T}$ defined by the vanishing of $\ell=\ell_{\alpha_1}+\cdots+\ell_{\alpha_n}$.   The WP distance of a point $p$ to $\mathcal S$ satisfies in terms of $\ell(p)$: in general $d_{WP}(p,\mathcal S)\le (2\pi\ell)^{1/2}$ and 
locally for $\ell$ small, $d_{WP}(p,\mathcal S)=(2\pi\ell)^{1/2}+O(\ell^2)$.  
\end{theorem}

\begin{corollary}
\label{bounds}
For $\beta$ represented by a simple free homotopy class the WP gradient of $\ell_\beta$ satisfies $\bigl<grad\,\ell_\beta,grad\,\ell_\beta\bigr>_{WP} \ge \frac{2}{\pi}\ell_\beta$.
As above, for $\alpha_1,\dots,\alpha_n$ and $\beta$ represented by disjoint distinct free homotopy classes: for $\gamma(s),\, 0\le s\le s_0$ the unit-speed distance-realizing WP geodesic connecting $\mathcal S$ to $p$ the derivatives of $(2\pi\ell)^{1/2}$ and $\ell_{\beta}$ along $\gamma$ satisfy $\frac{d}{ds}(2\pi\ell)^{1/2}(\gamma(s))\ge 1$ and $\frac{d}{ds}\ell_{\beta}(\gamma(s))\ge 0$.
\end{corollary}

G. Riera has recently obtained an exact formula for $\bigl<grad\,\ell_\alpha,grad\,\ell_\beta\bigr>_{WP}$ as an infinite sum for the lengths of the minimal geodesics connecting $\alpha$ to $\beta$ \cite{Rier}.   The above lower bound for $\bigl<grad\,\ell_\beta,grad\,\ell_\beta\bigr>_{WP}$ also follows from his formula.  The lower bound and the bound (\ref{gradbd}) can be combined to show that the {\em injectivity radius} $inj_{WP}$ (the minimal distance to a proper sub stratum in $\overline{\mathcal T}$) of $\mathcal T$ is comparable to the square root of the least geodesic-length.  In particular the bounds provide for positive constants $c_*,c_{**}$ and $c_{***}$ such that for $\ell=\ell_{\alpha_1}+\dots+\ell_{\alpha_n}$, $\ell\le c_*$ then $c_*\ell\le\bigl<grad\,\ell,grad\,\ell\bigr>_{WP}\le c_{**}\ell$ and for $\ell(R)\ge c_*$, $d_{WP}(R,R')<c_{***}$ then $\ell(R')\ge c_*/2$.  The overall bound $c'\,inj_{WP}\le (sys_{rel})^{1/2}\le c''\,inj_{WP}$ for positive constants is a consequence of the fact that $inj_{WP}$ and $(sys_{rel})^{1/2}$ are comparable for small values and are bounded in general.

We now present in detail two further applications for the behavior of WP geodesics.  The first is {\em Brock's approximation by rays to maximally noded surfaces} \cite{Brkwpvs}.  J. Brock noted that the $CAT(0)$ geometry and the observation of Bers on bounded partitions \cite{Bersdeg} provide for an approximation of infinite WP geodesics.  First note that the (incomplete) finite length WP geodesics from a point of $\mathcal T$ to the marked noded Riemann surfaces can be extended to include their endpoints in $\overline{\mathcal T}$.  As a consequence of the $CAT(0)$ geometry 
the initial unit tangents for such geodesics from a point to a stratum provide for a Lipschitz map from the stratum to the unit tangent sphere of the point.  Accordingly the image of $\overline{\mathcal T}-\mathcal T$ in
each unit tangent sphere has measure zero and consequently the infinite length geodesic rays have tangents dense in each tangent sphere.  In particular to approximate rays it suffices to approximate the infinite length rays.

From the result of Bers there is a positive constant $L_{g,n}$ depending only on the genus and number of punctures such 
that each surface has a maximal collection of simple closed curves $\alpha_1,\dots,\alpha_{3g-3+n}$ (a partition) with total geodesic-length bounded by $L_{g,n}$.  By Corollary \ref{bounds} each point of $\mathcal T$ is at most distance  $(2\pi L_{g,n})^{1/2}$ to a maximally noded Riemann surface.  To approximate an infinite ray $\gamma$ in $\mathcal T$ with initial point $p$, consider a point $q$ on the ray with $d_{WP}(p,q)$ large.  The point $q$ is at distance at most $(2\pi L_{g,n})^{1/2}$ to a maximally noded Riemann surface $q^\#$.  Since $\overline{\mathcal T}$ is a unique geodesic space the triple $(p,q,q^\#)$ determines a geodesic triangle.  The comparison to a Euclidean triangle provides that the initial angle between $\widehat{pq}$ and $\widehat{pq^\#}$ is $O(L_{g,n}^{1/2}\,d_{WP}(p,q)^{-1})$ \cite[Chap. II.1, Prop. 1.7]{BH}.  Further since $\gamma$ has infinite length the geodesic differential equation on $\mathcal T$ ensures that {\em close initial tangents provides for close initial segments.}  It further follows that initial segments of $\widehat{pq}$ and $\widehat{pq^\#}$ are close.  The considerations are summarized with the following.
\begin{theorem}
In $\mathcal T$ the infinite length geodesic rays and the rays to maximally noded Riemann surfaces each have initial tangents dense in each tangent space.  
\end{theorem}
Brock discovered that the situation for finite rays is different:  convergence of initial ray segments to finite rays does not provide for convergence of entire rays \cite{Brkwpvs}, \cite[Sec. 7]{Wlcomp}.  Rays approximating a finite ray can behave in a special way.  At this time an additional question is to understand infinite rays asymptotic to a stratum.  

As our second application we present a construction for asymptotic rays.  
Begin with the Teichm\"{u}ller space $\mathcal T'$ for a surface with $n>0$ punctures and $\mathcal A$ the axis for a pseudo Anosov mapping class.  The existence and uniqueness of a pseudo Anosov axis was first established in the work of G. Daskalopoulos and R. Wentworth \cite[Thm. 1.1]{DW2}.  In \cite[Thm. 25]{Wlcomp} the result was also obtained as an application of the classification of limits of geodesics and the general study of {\em translation length} \cite[Chap. II.6]{BH}.  Let $\{R'\}$ be the family of marked Riemann surfaces forming the axis $\mathcal A$ and $R''$ a particular Riemann  surface with $n$ punctures.  We view $\{R'\}$ as a surface bundle over $\mathcal A$ and $R''$ as a bundle over a point.  
We introduce a formal bijective pairing of the punctures of $\{R'\}$ with the punctures of $R''$ and consider the sum of surface bundles along fibers $\{R'\}+R''$ as a family of marked noded Riemann surfaces.  The nodes are the paired punctures.  For $g$ the formal genus of the family let $\mathcal T$ be the Teichm\"{u}ller space of genus $g$ surfaces with the length function $\ell=\ell_{\alpha_1}+\dots+\ell_{\alpha_n}$ defining the stratum $\mathcal S$ containing $\{R'\}+R''$ (the nodes have free homotopy classes $\alpha_1,\dots,\alpha_n$).  Further let $\gamma$ (reducible and partially pseudo Anosov) be an element of the mapping class group for $\mathcal T$ given as a sum of the pseudo Anosov (for $\{R'\}$) and the identity (for $R''$).  The mapping class $\gamma$ fixes $\alpha_1,\dots,\alpha_n$ and the action of $\gamma$ extends to $\mathcal S$ with the extension acting as the product of the pseudo Anosov on $\mathcal T'$ and the identity on $\mathcal T(R'')$.

We proceed and describe the construction of a geodesic ray in $\mathcal T$ asymptotic to $\mathcal S$.  First observe that the relative systole is periodic along $\mathcal A$ and consequently that $sys_{rel}$ is bounded below along $\mathcal A$ by a positive constant $c$.
It now follows from the gradient bound (\ref{gradbd}) that there exists a positive constant $\delta$ such that any surface $R$ of $\mathcal T$ closer in $\overline{\mathcal T}$ to $\{R'\}+R''$ than $\delta$ satisfies $0\le\ell_{\alpha_j}<c/3$ and for $\beta\ne \alpha_1,\dots,\alpha_n$ (or a power of an $\alpha_j$) then $\ell_{\beta}(R)\ge 2c/3$.  In particular the only {\em short} primitive geodesics on such an $R$ are $\alpha_1,\dots,\alpha_n$.  We are ready to form the candidate ray asymptotic to $\mathcal S$ by a limiting process.  For a sequence of points along $\{R'\}+R''$ tending to forward infinity, connect the reference point $R$ by a WP geodesic to each point of the sequence.  The point $R$ in 
$\mathcal T$ has relatively compact neighborhoods and consequently we can select a convergent subsequence of the connecting geodesics.  Denote the resulting limit as $\mathcal G$.  We will verify that the limit is an infinite ray.  We are interested in the behavior of  three functions on $\mathcal G$: $\ell_{\alpha_j}$, $d_{WP}(\ ,\{R'\}+R'')$ and $d_{WP}(\ ,\mathcal S)$.  On each geodesic connecting $R$ to a point of $\{R'\}+R''$ each of the functions is convex (see the above on orthogonal projections). Further each function vanishes at the far endpoint of each connecting geodesic.  It follows that each function is strictly decreasing on each connecting geodesic and consequently that each function is non increasing on the limit $\mathcal G$. Now the classification of geodesic (with lengths tending to infinity) limits provides that either: a limit is an infinite ray or at a fixed distance from the basepoint: the limiting rays successively approach and then strictly recede from a stratum \cite[Prop. 23]{Wlcomp}.  As already noted along $\mathcal G$ the only possible {\em small} geodesic-lengths have non increasing length functions: the second limiting behavior is precluded and consequently the limit is an infinite ray. 

We will now show that $\ell_{\alpha_j}$ and $d_{WP}(\ ,\mathcal S)$ tend to zero along $\mathcal G$.  For this sake consider $\mathcal N_{\delta}$: the points in $\overline{\mathcal T}$ at distance at most $\delta$ from $\{R'\}+R''$.  The closed set $\mathcal N_{\delta}$ is stabilized by the action of the mapping class $\gamma$, as well as by the Dehn twists $\tau_j$ about $\alpha_j$.  Since the only possible short primitive geodesicsf for a surface in $\mathcal N_{\delta}$ are $\alpha_1,\dots,\alpha_n$ it can be shown that the quotient of $\mathcal N_{\delta}$ by the action of the group generated by $\gamma$ and the $\tau_j$ is compact.  Now for a sequence of points along $\mathcal G$ tending to infinity consider the associated sequence of forward direction rays.  
Since the quotient $\mathcal N_{\delta}$ is compact we can select a convergent subsequence of rays translated by appropriate compositions with powers of $\gamma$ and the $\tau_j$ \cite[Prop. 23]{Wlcomp}.  The resulting limit is a geodesic $\mathcal G_0$ in $\mathcal N_{\delta}$.  Since $\ell_{\alpha_j}$, $d_{WP}(\ ,\{R'\}+R'')$ and $d_{WP}(\ ,\mathcal S)$ are non increasing along $\mathcal G$, each function has a limit along $\mathcal G$ and consequently each function is actually constant on $\mathcal G_0$.  In particular each $\ell_{\alpha_j}$ is constant on $\mathcal G_0$; Theorem \ref{convex} provides that each $\ell_{\alpha_j}$ vanishes on $\mathcal G_0$.  It further follows from Theorem \ref{dist} that $d_{WP}(\ ,\mathcal S)$ vanishes on $\mathcal G_0$ and thus that $\mathcal G$ is asymptotic to $\mathcal S$, as proposed.  Finally in closing we note that if the Teichm\"{u}ller space $\mathcal T(R'')$ is not a singleton then the product $\mathcal T(R')\times\mathcal T(R'')$ contains {\em Euclidean flats} and $\gamma$ stabilizes parallel lines $\{R'\}+R'',\,\{R'\}+R'''$.  We expect families of asymptotic rays in this case.    

\bibliographystyle{alpha}

\clearpage

\bigskip

\end{document}